# On the Existence of Lagrange Multipliers in Distribution Network Reconfiguration Problems

Rong-Peng Liu, *Member, IEEE*, Yue Song, *Member, IEEE*, Xiaozhe Wang, *Senior Member, IEEE*, and Bo Zeng, *Member, IEEE*

*Abstract*—Distribution network reconfiguration (DNR) is an effective approach for optimizing distribution network operation. However, the DNR problem is computationally challenging due to the mixed-integer non-convex nature. One feasible approach for addressing this challenge is to reformulate binary line on/off state variables as (continuous) non-convex constraints, leading to a nonconvex program. Unfortunately, it remains unclear whether this formulation satisfies the Karush-Kuhn-Tucker (KKT) conditions at locally optimal solutions. In this brief, we study the existence of Lagrange multipliers in DNR problems and prove that under mild assumptions, Lagrange multipliers exist for the DNR model at every locally optimal solution almost surely such that the KKT conditions hold.

## I. Introduction

DISTRIBUTION networks are critical infrastructures for delivering electricity from power transmission networks to end-users. Distribution network reconfiguration (DNR) [1] is an effective approach to enhance system reliability against fault [1] and cyber-physical attacks [2] by isolating faults, preventing potential cascading failures [3]. In addition, DNR is also effective in enhancing network performance, e.g., loss reduction [4], voltage regulation [5], congestion mitigation [6], and stability enhancement [7]. Physically, DNR is to adjust the topology of a distribution network by altering the on/off states of switches mounted in feeders [8]. After reconfiguration, the distribution network should maintain a radial structure. In practice, DNR can be actualized by manipulating remote-controlled switches [8], which is readily implementable and cost-effective for optimizing network states.

This work was supported in part by the Fonds de recherche du Québec-secteur Nature et technologies (FRQ-Secteur NT) under Grant FRQ-Secteur NT 334636, Grant FRQ-Secteur NT 367013, Grant FRQ-Secteur NT PR-298827, and Grant FRQ-Secteur NT 2023-NOVA-314338, in part by the National Natural Science Foundation of China under Grant 62088101, and in part by the Shanghai Municipal Science and Technology Major Project under Grant 2021SHZDZX0100.

R. Liu and X. Wang are with the Department of Electrical and Computer Engineering, McGill University, Montreal, QC H3A 0E9, Canada (e-mail: rpliu@eee.hku.hk/rongpeng.liu@mail.mcgill.ca; xiaozhe.wang2@mcgill.ca).

Y. Song is with the Department of Control Science and Engineering, Tongji University, Shanghai 201804, China, also with the National Key Laboratory of Autonomous Intelligent Unmanned Systems (Tongji University), Shanghai 201210, China, and also with the Shanghai Institute of Intelligent Science and Technology, Tongji University, Shanghai 201210, China (e-mail: ysong@tongji.edu.cn).

B. Zeng is with the Department of Industrial Engineering and the Department of Electrical and Computer Engineering, University of Pittsburgh, Pittsburgh, PA 15106 USA (e-mail: bzeng@pitt.edu).

Despite these advantages, DNR problems, mixed-integer nonconvex programs (MINCPs), are computationally intractable due to nonconvex power flow constraints and binary line (on/ off) state constraints. One method to address the (computational) challenge is to simplify the power flow constraints, e.g., linear approximation [1], [5] and convex relaxation [4], [6], [9], leading to mixed-integer convex programs. Another method is to employ machine learning-based approaches to solve DNR problems directly [7]. However, the solutions derived by either (simplification- or machine learning-based) method may not be feasible. As a critical energy infrastructure, the feasibility of a reconfiguration strategy, i.e., a solution, is of great significance in practice and should be put higher priority. Differently, reference [10] tries to address the challenge by reformulating binary line state variables as (continuous) nonconvex constraints [11], resulting in nonconvex programs (NCPs). This method ensures both the feasibility and (local) optimality of a solution and is considered reasonable in coping with practical DNR problems. In this work, we focus on this (continuous) nonconvex DNR formulation in light of its practicability and, hereinafter, denote it as *the DNR model*.

The Karush-Kuhn-Tucker (KKT) conditions [12] play a critical role in addressing the nonconvex DNR model. First, it serves as the foundation of some algorithms for solving nonconvex optimization problems, such as the primal-dual interior point method [13]. Second, the KKT conditions are effective in analyzing the optimality of a solution [12]. Specifically, there exist Lagrange multipliers such that the optimal solution of an optimization problem must satisfy some constraints (presented in Section II), which is referred to as the KKT conditions, if at least one of the constraint qualifications (CQs) is satisfied [14]. Note that the CQs are crucial for ensuring the existence of Lagrange multipliers such that the KKT conditions hold. If none of the CQs is satisfied, an optimal solution cannot meet the KKT conditions (due to the non-existence of Lagrange multipliers). As a result, KKT-based solution and analysis methods are inapplicable. Unfortunately, the CQs are hard to verify, especially for nonconvex optimization problems, and are usually ignored [15] or considered as an assumption [16] in some previous works.

In this brief, we study the existence of Lagrange multipliers in the DNR model. For the first time, we prove that the linear independent constraint qualification (LICQ), one of the most frequently employed CQs for nonconvex optimization problems, does not hold for the DNR model. Since the LICQ is the weakest CQ that ensures the uniqueness (and the existence) of Lagrange multipliers [14], unique Lagrange multipliers do not exist for the DNR model. Note that although not unique, (multiple) Lagrange multipliers may still exist, i.e., some weaker CQs may be satisfied. Then, for the first time, we propose an auxiliary model and, based on it, prove that under mild assumptions, Lagrange multipliers, although not unique, exist for the DNR model at every locally optimal solution almost surely such that the KKT conditions hold. This conclusion renders a solid theoretical foundation for the application of KKT-based solution methods, e.g., the primal-dual interior point method, for solving the DNR model and the analysis of the solution optimality via the KKT conditions.

The remainder of this brief is organized as follows. Section II introduces the DNR model and the KKT conditions. Section III studies the existence of Lagrange multipliers in the DNR model. Section 0 conducts case studies. Section V concludes this brief and discusses future works.

## II. PROBLEM FORMULATION

This section presents the DNR formulation. Before proceeding, we introduce the scope of this work.

1) This work studies three-phase balanced distribution networks and thus focuses on a single phase. Three-phase unbalanced distribution networks are beyond the scope of this work.

2) This work studies steady-state DNR and does not consider the impact of system dynamics caused by network reconfiguration.

3) This work studies deterministic DNR and does not consider uncertainties, e.g., power load uncertainty. DNR under uncertainties remains our future work.

TABLE I
NOMENCLATURE

| Notation | Explanation |
|---|---|
| $\mathcal{L}_n^+/\mathcal{L}_n^-$ | Set of (power transmission) lines connecting to bus $n$ with a presupposed inflow/outflow direction to/from bus $n$. [Set] |
| $\mathcal{N}/\mathcal{L}$ | Set of buses/lines. [Set] |
| $\mathcal{N}_n^+$ | Set of neighboring buses that directly connect to bus $n$ with a presupposed inflow direction to bus $n$. [Set] |
| $C_{mn}$ | Cost of line losses in line $mn$. [Parameter] |
| $d_n^P/d_n^Q$ | Real/reactive load at bus $n$. [Parameter] |
| $P_{mn}^{max}/Q_{mn}^{max}$ | Limit of real/reactive power flow in line $mn$. [Parameter] |
| $r_{mn}/x_{mn}$ | Resistance/reactance of line $mn$. [Parameter] |
| $V_n^{min}/V_n^{max}$ | Limit of voltage magnitude at bus $n$. [Parameter] |
| $l_{mn} \in \mathbb{R}$ | Squared current magnitude in line $mn$. [Variable] |
| $p_{mn}/q_{mn} \in \mathbb{R}$ | Real/reactive power flow in line $mn$. [Variable] |
| $p_n/q_n \in \mathbb{R}$ | Real/reactive power injection at bus $n$. [Variable] |
| $u_{mn} \in \mathbb{R}$ | State of line $mn$. 1 if line $mn$ is switched on, and 0 otherwise. [Variable] |
| $v_n \in \mathbb{R}$ | Squared voltage magnitude at bus $n$. [Variable] |

### A. Preliminaries

*Notations*: In this brief, we adopt i) regular uppercase or lowercase letters to denote scalars, ii) bold uppercase and lowercase letters to denote matrices and column vectors, respectively, iii) calligraphic or hollow uppercase letters to denote sets. Vectors are all column vectors, and superscript "T" denotes the transpose manipulation. Function $col(\cdot)$ denotes the reshape mapping of multiple vectors or matrices. For example, vector $col(\mathbf{a}, b, \mathbf{c}) = [\mathbf{a}^T\ b\ \mathbf{c}^T]^T$ is a vector consisting of vectors $\mathbf{a}$, b, and $\mathbf{c}$, and matrix $col(\mathbf{A}, \mathbf{B}) = [\mathbf{A}^T\ \mathbf{B}^T]^T$ is a matrix consisting of matrices $\mathbf{A}$ and $\mathbf{B}$, where matrices $\mathbf{A}$ and $\mathbf{B}$ should have the same number of columns. Function $|\cdot|$ denotes the number of elements in a set.

*Nomenclature*: We list the main sets, parameters, and variables, as well as their explanations, used in this brief in Table I.

### B. Distribution Network Reconfiguration (DNR) Model

The DNR model is presented as follows.

$$\min_{\mathbf{z}} \sum_{ij \in \mathcal{L}} r_{ij} l_{ij} \tag{1a}$$

$$\text{s.t.}\ p_j = \sum_{ij \in \mathcal{L}_j^+} (r_{ij}l_{ij} - p_{ij}) + \sum_{jk \in \mathcal{L}_j^-} p_{jk},\ \forall j \in \mathcal{N}, \tag{1b}$$

$$q_j = \sum_{ij \in \mathcal{L}_j^+} (x_{ij}l_{ij} - q_{ij}) + \sum_{jk \in \mathcal{L}_j^-} q_{jk},\ \forall j \in \mathcal{N}, \tag{1c}$$

$$v_i - v_j \geq 2(r_{ij}p_{ij} + x_{ij}q_{ij}) - (r_{ij}^2 + x_{ij}^2)l_{ij} - M(1 - u_{ij}),\ \forall ij \in \mathcal{L}, \tag{1d}$$

$$v_i - v_j \leq 2(r_{ij}p_{ij} + x_{ij}q_{ij}) - (r_{ij}^2 + x_{ij}^2)l_{ij} + M(1 - u_{ij}),\ \forall ij \in \mathcal{L}, \tag{1e}$$

$$v_i l_{ij} = p_{ij}^2 + q_{ij}^2,\ \forall ij \in \mathcal{L}, \tag{1f}$$

$$-u_{ij}P_{ij}^{max} \leq p_{ij} \leq u_{ij}P_{ij}^{max},\ \forall ij \in \mathcal{L}, \tag{1g}$$

$$-u_{ij}Q_{ij}^{max} \leq q_{ij} \leq u_{ij}Q_{ij}^{max},\ \forall ij \in \mathcal{L}, \tag{1h}$$

$$(V_i^{min})^2 \leq v_i \leq (V_i^{max})^2,\ \forall i \in \mathcal{N}, \tag{1i}$$

$$u_{ij}(u_{ij} - 1) = 0,\ \forall ij \in \mathcal{L}, \tag{1j}$$

$$\sum_{ij \in \mathcal{L}_j^+} u_{ij} + \sum_{jk \in \mathcal{L}_j^-} u_{jk} \geq 1,\ \forall j \in \mathcal{N}, \tag{1k}$$

$$\sum_{j \in \mathcal{N}} \sum_{ij \in \mathcal{L}_j^+} u_{ij} = |\mathcal{N}| - 1. \tag{1l}$$

Variable $\mathbf{z} \in \mathbb{R}^{3|\mathcal{N}|+4|\mathcal{L}|} = col(p_j, q_j, p_{ij}, q_{ij}, l_{ij}, v_j, u_{ij})$, $\forall ij \in \mathcal{L}$, $j \in \mathcal{N}$. Objective function (1a) minimizes aggregated line losses. Constraints (1b)-(1f) are a modified DistFlow model based on the original DistFlow model [1] considering the impact of line states (i.e., switched on or off). M is a large constant and could invalidate constraints (1d) and (1e) when corresponding lines are switched off. Constraints (1g) and (1h) regulate real and reactive power flow, respectively, considering the impact of line states. Constraint (1i) restricts squared voltage magnitudes. Constraint (1j) confines the states of each line to be either 0 (switched off) or 1 (switched on). Constraints (1k) and (1l) ensure the fully connected radial structure of a distribution network after reconfiguration [5], i.e., $|\mathcal{N}|$ - 1 lines are switched-on and the rest $|\mathcal{L}|$ - $|\mathcal{N}|$ + 1

lines (if applicable) are switched-off.

DNR model (1) is a nonconvex optimization problem due to DistFlow constraint (1f) and line state constraint (1j) and is computationally challenging. The interior point method is effective in solving this model, where the KKT conditions play a critical role. To introduce the KKT conditions, we first present the compact form for (1).

$$\mathbf{P}: \min_{\mathbf{z}} f(\mathbf{z}) \quad (2a)$$
$$\text{s.t.} \quad \mathbf{g}(\mathbf{z}) = \mathbf{0}, \quad (\boldsymbol{\lambda}) \quad (2b)$$
$$\mathbf{h}(\mathbf{z}) \leq \mathbf{0}, \quad (\boldsymbol{\mu}) \quad (2c)$$
$$\mathbf{z} \in \mathbb{R}^Z. \quad (2d)$$

Objective function $f(\cdot)$ denotes the objective function in (1a), equality constraints $\mathbf{g}(\cdot): \mathbb{R}^Z \to \mathbb{R}^G$ consist of (1b), (1c), (1f), (1j), (1l), and inequality constraint $\mathbf{h}(\cdot): \mathbb{R}^Z \to \mathbb{R}^H$ consists of (1d), (1e), (1g)-(1i), (1k). $Z = 3|\mathcal{L}| + 3|\mathcal{N}|$. Vectors $\boldsymbol{\lambda}$ and $\boldsymbol{\mu}$ are Lagrange multipliers (vectors) for equality and inequality constraints, respectively. For locally optimal solution $\mathbf{z}^*$, we can find Lagrange multipliers $\boldsymbol{\lambda}^*$ and $\boldsymbol{\mu}^*$ such that the following conditions, i.e., the *KKT conditions*,

$$-\nabla f(\mathbf{z}^*)^T = (\boldsymbol{\lambda}^*)^T \nabla \mathbf{g}(\mathbf{z}^*) + (\boldsymbol{\mu}^*)^T \nabla \mathbf{h}(\mathbf{z}^*), \quad (3a)$$
$$\mathbf{g}(\mathbf{z}^*) = \mathbf{0}, \quad (3b)$$
$$\mathbf{h}(\mathbf{z}^*) \leq \mathbf{0}, \quad (3c)$$
$$\boldsymbol{\mu}^* \geq \mathbf{0}, \quad (3d)$$
$$(\boldsymbol{\mu}^*)^T \mathbf{h}(\mathbf{z}^*) = 0 \quad (3e)$$

hold if constraints (2b) and (2c) satisfy one of the CQs [14] at $\mathbf{z}^*$. In other words, Lagrange multipliers that satisfy the KKT conditions may not exist even at an optimal solution due to the violation of CQs. In the next section, we study if Lagrange multipliers $\boldsymbol{\lambda}^*$ and $\boldsymbol{\mu}^*$ exist for DNR model (1).

## III. EXISTENCE OF LAGRANGE MULTIPLIERS

In this section, we analyze the existence of Lagrange multipliers in DNR model (1). As aforementioned, for the optimal solution of a nonconvex optimization problem, Lagrange multipliers exist such that the KKT conditions are satisfied if the constraints at the optimal solution meet one of the CQs [14]. Among them, the LICQ has been proved to be the weakest CQ that can ensure both the existence and the uniqueness of Lagrange multipliers [14] and thus is broadly adopted as an assumption in optimization problems. Specifically, the LICQ requires that the gradients of all equality constraints and active inequality constraints at a feasible solution, including an optimal solution, are linearly independent. For example, DNR model (1) satisfies the LICQ if we have

$$rank(col(\nabla \mathbf{g}(\mathbf{z}), \nabla \mathbf{h}_{\mathcal{A}}(\mathbf{z}))) = G + |\mathcal{A}|, \quad (4)$$

where $\mathbf{h}_{\mathcal{A}}(\cdot)$ denotes active inequality constraints at $\mathbf{z}$, and $|\mathcal{A}|$ is the number of active inequality constraints. Unfortunately, DNR model (1) does not satisfy the LICQ. Specifically, we have the following lemma.

*Lemma 1:* The LICQ does not hold for the constraints in (1) at any feasible solution.

*Proof:* Without loss of generality, we assume DNR model (1) has $|\mathcal{N}|$ buses and $|\mathcal{L}|$ lines (including both switched-on and switched-off lines), where $|\mathcal{L}| \geq |\mathcal{N}| - 1$. The partial derivative of constraints (1j) and (1l) with respect to variables $\mathbf{u}_{|\mathcal{L}|}$ is presented in the following matrix form:

$$\left[ Diag(2\mathbf{u}_{|\mathcal{L}|} - \mathbf{1}_{|\mathcal{L}|}) \mid \mathbf{1}_{|\mathcal{L}|} \right]^T. \quad (5)$$

Matrix $Diag(2\mathbf{u}_{|\mathcal{L}|} - \mathbf{1}_{|\mathcal{L}|})$ is a diagonal matrix, and its diagonal elements consist of vector $2\mathbf{u}_{|\mathcal{L}|} - \mathbf{1}_{|\mathcal{L}|}$. Vector $\mathbf{u}_{|\mathcal{L}|}$ is an $|\mathcal{L}|$ dimensional vector. Note that the value of each element in $\mathbf{u}$ can only be 0 or 1. Vector $\mathbf{1}_{|\mathcal{L}|}$ is an all-one vector. For any feasible solution of (1), matrix (5) has a general form, i.e.,

$$\begin{bmatrix} Diag(\mathbf{1}_{|\mathcal{N}|-1}) & \mathbf{0}_{(|\mathcal{N}|-1)*(|\mathcal{L}|-|\mathcal{N}|+1)} \\ \mathbf{0}^T_{(|\mathcal{N}|-1)*(|\mathcal{L}|-|\mathcal{N}|+1)} & Diag(-\mathbf{1}_{|\mathcal{L}|-|\mathcal{N}|+1}) \end{bmatrix} \mathbf{1}_{|\mathcal{L}|} \Bigg]^T. \quad (6)$$

This is because a distribution network should keep a fully connected radial structure after reconfiguration, i.e., $|\mathcal{N}| - 1$ lines are switched-on and the rest $|\mathcal{L}| - |\mathcal{N}| + 1$ lines (if applicable) are switched-off. It is easy to validate that the LICQ is not satisfied for the constraints in (1) at any feasible solution, as the rank of matrix (6) is always smaller than $|\mathcal{L}| + 1$, i.e., the number of rows of matrix (6), indicating the gradients of constraints (1j) and (1l) at any feasible solution of (1) are linearly dependent. This completes the proof. ∎

Based on the above proof, we immediately derive the following corollary for another important CQ, i.e.,

*Corollary 1*: The Mangasarian-Fromovitz constraint qualification (MFCQ) does not hold for the constraints in (1) at any feasible solution.

*Remark 1*: The MFCQ, weaker than the LICQ, is another widely used CQ. Lemma 1 and Corollary 1 demonstrate that either the LICQ or the MFCQ is not satisfied for DNR model (1) due to reconfiguration constraints (1j)-(1l). Note that the violation of these two CQs does not indicate that Lagrange multipliers do not exist at an optimal solution. For example, reference [14] indicates that the LICQ is the weakest CQ that ensures their *existence* and *uniqueness*. Furthermore, reference [13] points out that there are weaker CQs that only guarantee the existence of Lagrange multipliers at an optimal solution. Thus, the existence of Lagrange multipliers in (1) should be further analyzed.

Before proceeding, we introduce two models, $\mathbf{P}(\boldsymbol{\delta})$ and $\mathbf{P}_{R\_E}(\mathbf{d}, \mathcal{T})$. Fig. 1 presents the relations between DNR model (1), $\mathbf{P}(\boldsymbol{\delta})$ and $\mathbf{P}_{R\_E}(\mathbf{d}, \mathcal{T})$.

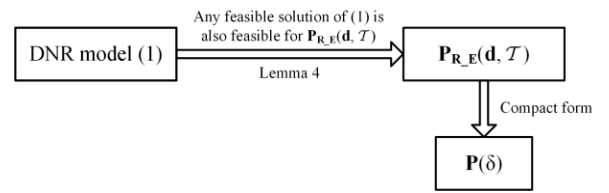

Fig. 1. Relations between DNR model (1), $\mathbf{P}(\boldsymbol{\delta})$, and $\mathbf{P}_{R\_E}(\mathbf{d}, \mathcal{T})$.

First, we introduce general optimization model $\mathbf{P}(\boldsymbol{\delta})$, where

$$\mathbf{P}(\boldsymbol{\delta}): \min_{\mathbf{z}} w(\mathbf{z}) \qquad (7a)$$
$$\text{s.t. } \mathbf{g}_1(\mathbf{z}) = \mathbf{0},\ \mathbf{g}_2(\mathbf{z}\mid\boldsymbol{\delta}) = \mathbf{0}, \quad (\boldsymbol{\iota}) \qquad (7b)$$
$$\mathbf{h}_1(\mathbf{z}) \leq \mathbf{0},\ \mathbf{h}_2(\mathbf{z}\mid\boldsymbol{\delta}) \leq \mathbf{0}, \quad (\boldsymbol{\kappa}) \qquad (7c)$$
$$\mathbf{z} \in \mathbb{R}^Z. \qquad (7d)$$

Hyperparameter $\boldsymbol{\delta} \in \mathcal{W}$, where set $\mathcal{W}$ is a nonempty open set with non-zero measure. $w(\cdot)$ is the objective function. $\mathbf{g}_1(\cdot)$: $\mathbb{R}^Z \to \mathbb{R}^{G_1}$ and $\mathbf{h}_1(\cdot)$: $\mathbb{R}^Z \to \mathbb{R}^{H_1}$, termed *fixed constraints* [13], are equality and inequality constraints, respectively, without hyperparameter $\boldsymbol{\delta}$. $\mathbf{g}_2(\cdot)$: $\mathbb{R}^Z \times \mathcal{W} \to \mathbb{R}^{G_2}$ and $\mathbf{h}_2(\cdot)$: $\mathbb{R}^Z \times \mathcal{W} \to \mathbb{R}^{H_2}$, termed *variable constraints* [13], are equality and inequality constraints, respectively, with hyperparameter $\boldsymbol{\delta}$. Operator "$\times$" denotes the Cartesian product. Set $\mathcal{Z} := \{\mathbf{z}\mid \mathbf{g}_1(\mathbf{z}) = \mathbf{0}, \mathbf{h}_1(\mathbf{z}) \leq \mathbf{0}\}$. Variables $\boldsymbol{\iota}$ and $\boldsymbol{\kappa}$ are Lagrange multipliers for equality constraints $col(\mathbf{g}_1(\mathbf{z}) = \mathbf{0}, \mathbf{g}_2(\mathbf{z}\mid\boldsymbol{\delta}) = \mathbf{0})$ and inequality constraints $col(\mathbf{h}_1(\mathbf{z}) \leq \mathbf{0}, \mathbf{h}_2(\mathbf{z}\mid\boldsymbol{\delta}) \leq \mathbf{0})$, respectively.

For $\mathbf{P}(\boldsymbol{\delta})$, we have an important theorem.

*Lemma 2 (Theorem 2 [13]):* Assume that: i) objective function $w(\cdot)$ and variable constraints $\mathbf{g}_2(\cdot)$ and $\mathbf{h}_2(\cdot)$ are (1-time) continuously differentiable; ii) fixed constraints $\mathbf{g}_1(\mathbf{z})$ and $\mathbf{h}_1(\mathbf{z})$ are $c$-times continuously differentiable for every $\mathbf{z} \in \mathcal{Z}$, where $c > \max(0, Z - G_1 - G_2)$; iii) the LICQ holds for every $\mathbf{z} \in \mathcal{Z}$ with respect to $\mathbf{g}_1(\mathbf{z})$ and active $\mathbf{h}_1(\mathbf{z})$ [13]; iv) for every $\mathbf{z} \in \mathcal{Z}$ and every $\boldsymbol{\delta} \in \mathcal{W}$, the map $\boldsymbol{\delta} \to col(\mathbf{g}_2(\mathbf{z}\mid\boldsymbol{\delta}), \mathbf{h}_2(\mathbf{z}\mid\boldsymbol{\delta}))$ has rank $G_2 + H_2$. Then, for almost every $\boldsymbol{\delta} \in \mathcal{W}$, the LICQ holds for the constraints in $\mathbf{P}(\boldsymbol{\delta})$ at every feasible solution.

Please refer to [13] for the detailed proof of Lemma 2. Then, we introduce an auxiliary model $\mathbf{P}_{R\_E}(\mathbf{d}, \mathcal{T})$, and its formulation as well as detailed explanations of hyperparameters $\mathbf{d}$ and $\mathcal{T}$ is presented in Appendix.A. For any fixed $\mathcal{T} \in \mathbb{T}$, we use $\mathbf{P}(\boldsymbol{\delta})$ to denote the compact form of $\mathbf{P}_{R\_E}(\mathbf{d}, \mathcal{T})$. Based on Lemma 2, we have

*Lemma 3:* Assume that for every $\mathcal{T} \in \mathbb{T}$ in $\mathbf{P}_{R\_E}(\mathbf{d}, \mathcal{T})$, the LICQ holds at every $\mathbf{z} \in \mathcal{Z}$ with respect to fixed constraints $\mathbf{g}_1(\cdot)$ and active $\mathbf{h}_1(\cdot)$. Then, for every $\mathcal{T} \in \mathbb{T}$ and almost every $\mathbf{d} \in \mathcal{D}$, the LICQ holds for the constraints in $\mathbf{P}_{R\_E}(\mathbf{d}, \mathcal{T})$ at every feasible solution.

*Proof:* As aforementioned, for any fixed $\mathcal{T} \in \mathbb{T}$, we use $\mathbf{P}(\boldsymbol{\delta})$ to denote the compact form of $\mathbf{P}_{R\_E}(\mathbf{d}, \mathcal{T})$, where hyperparameter $\boldsymbol{\delta}$ denotes $\mathbf{d}$, set $\mathcal{W}$ denotes $\mathcal{D}$, objective function $w(\cdot)$ denotes the objective function in (A.1), equality constraints $\mathbf{g}_1(\cdot)$ and $\mathbf{g}_2(\cdot)$ consist of (A.4)-(A.7) and (A.2)-(A.3), respectively, inequality constraint $\mathbf{h}_1(\cdot)$ consists of (A.8)-(A.10), and inequality constraint $\mathbf{h}_2(\cdot)$ does not exist. $\mathcal{Z} := \{\mathbf{z}\mid$ (A.4)-(A.10)$\}$. Since Lemma 3 assumes that assumption iii) in Lemma 2 is always satisfied, this proof is to verify assumptions i), ii), and iv) in Lemma 2. First, objective function $w(\cdot)$, i.e., (A.1), denotes any 1-time continuously differentiable function, and variable constraint $\mathbf{g}_2(\cdot)$, i.e., (A.2)-(A.3), is linear equality constraints, which are also 1-time continuously differentiable. Considering variable constraint $\mathbf{h}_2(\cdot)$ does not exist, assumption i) in Lemma 2 is satisfied. Second, for $\mathbf{P}_{R\_E}(\mathbf{d}, \mathcal{T})$, $G_1 = 3|\mathcal{L}| - |\mathcal{N}|$, $G_2 = 2|\mathcal{N}|$, and $Z = 3|\mathcal{L}| + 3|\mathcal{N}|$. Based on the second assumption in Lemma 2, parameter $c$ should be set to be larger than $2|\mathcal{N}|$. This is always true, since fixed constraints $\mathbf{g}_1(\cdot)$, i.e., (A.4)-(A.7), and $\mathbf{h}_1(\cdot)$, i.e., (A.8)-(A.10), consist of linear and quadratic equality constraints and affine inequality constraint, respectively. Last, for $\mathbf{P}_{R\_E}(\mathbf{d}, \mathcal{T})$, $G_2 = 2|\mathcal{N}|$ and $H_2 = 0$. The rank of map $\boldsymbol{\delta} \to col(\mathbf{g}_2(\mathbf{z}\mid\boldsymbol{\delta}), \mathbf{h}_2(\mathbf{z}\mid\boldsymbol{\delta}))$ is $2|\mathcal{N}|$, as we have

$$\frac{\partial \mathbf{g}_2(\mathbf{z}\mid\boldsymbol{\delta})}{\partial \boldsymbol{\delta}} = Diag(\mathbf{1}_{2|\mathcal{N}|}),$$

and inequality constraint $\mathbf{h}_2(\cdot)$ does not exist. This implies that assumption iv) in Lemma 2 is satisfied. Overall, all the assumptions in Lemma 2 are satisfied, and we derive Lemma 3. This completes the proof. ∎

*Remark 2:* Lemma 3 indicates that $\mathbf{P}_{R\_E}(\mathbf{d}, \mathcal{T})$ satisfies the LICQ almost surely. Note that $\mathbf{P}_{R\_E}(\mathbf{d}, \mathcal{T})$ does not allow reconfiguration manipulation and does not require a radial network topology. Next, we establish connections between $\mathbf{P}_{R\_E}(\mathbf{d}, \mathcal{T})$ and DNR model (1). Specifically, we have

*Lemma 4:* For every feasible solution of (1) with load $\hat{\mathbf{d}} \in \mathcal{D}$, i.e., $col(\hat{\mathbf{z}}, \hat{\mathbf{u}})$, there always exists a $\mathcal{T} \in \mathbb{T}$ such that this $\hat{\mathbf{z}}$ is also feasible for $\mathbf{P}_{R\_E}(\hat{\mathbf{d}}, \mathcal{T})$.

*Proof:* Without loss of generality, we assume that DNR model (1) and $\mathbf{P}_{R\_E}(\hat{\mathbf{d}}, \mathcal{T})$ have $|\mathcal{N}|$ buses and $|\mathcal{L}|$ lines, where $|\mathcal{L}| \geq |\mathcal{N}| - 1$. When $|\mathcal{L}| = |\mathcal{N}| - 1$, set $\mathbb{T}$ is a singleton. Its element corresponds to the (sole) reconfiguration decision $\mathbf{u}_{|\mathcal{L}|} = \mathbf{1}_{|\mathcal{L}|}$. For this case, it is easy to prove that Lemma 4 holds, since no reconfiguration manipulation is allowed. When $|\mathcal{L}| > |\mathcal{N}| - 1$, constraints (1b), (1c), (1f), and (1i) are the same as constraints, (A.2)-(A.3), (A.5), and (A.10), respectively. Every feasible $\hat{\mathbf{z}}$ for the former must be feasible for the latter. In addition, there always exists a $\mathcal{T} \in \mathbb{T}$ such that constraints (1g) and (1h) are the same as constraints (A.6)-(A.9). This is because set $\mathbb{T}$ consists of all connected radial networks. We can always find a $\mathcal{T} \in \mathbb{T}$ such that switched-on lines and switched-off lines in $\mathbf{P}_{R\_E}(\hat{\mathbf{d}}, \mathcal{T})$ are enforced by (A.8)-(A.9) and (A.6)-(A.7), respectively, i.e., the same as (1g) and (1h). Thus, every feasible $\hat{\mathbf{z}}$ for (1g) and (1h) must be feasible for (A.6)-(A.9). For constraints (1d) and (1e), they together become equality constraints for the $|\mathcal{N}| - 1$ switched-on lines that compose a connected radial network and are invalid for switched-off lines. In fact, constraint (A.4) only enforces the $|\mathcal{N}| - 1$ lines that constitute a connected radial network and does not restrict the other lines. Similarly, we can always find the same $\mathcal{T} \in \mathbb{T}$ (as the $\mathcal{T}$ in the analysis for (A.6)-(A.9)) such that the switched-off lines in (1) are not restricted by (A.4). This indicates that every feasible $\hat{\mathbf{z}}$ for (1d) and (1e) is also feasible for (A.4). Thus, Lemma 4 holds for the case when $|\mathcal{L}| > |\mathcal{N}| - 1$. This completes the proof. ∎

*Remark 3:* According to Lemma 1, it is the reconfiguration constraints (1j)-(1l) that induce the violation of the LICQ for DNR model (1). In addition, reconfiguration variable $\mathbf{u}$ also appears in constraints (1d), (1e), (1g), and (1h) and affects the feasible region of (1). As the first step for deriving $\mathbf{P}_{R\_E}(\hat{\mathbf{d}}, \mathcal{T})$, we *equivalently* replace reconfiguration constraints with hyperparameter $\mathcal{T} \in \mathbb{T}$. As the second step for deriving $\mathbf{P}_{R\_E}(\hat{\mathbf{d}}, \mathcal{T})$, we *equivalently* reformulate constraints (1g) and (1h), which contain reconfiguration decision $\mathbf{u}$, into (A.6)-(A.9) by means of hyperparameter $\mathcal{T} \in \mathbb{T}$ and set $\mathcal{L}(\mathcal{T})$. As the last step for deriving $\mathbf{P}_{R\_E}(\hat{\mathbf{d}}, \mathcal{T})$, we *equivalently* re-write constraints (1d) and (1e) as constraint (A.4) by means of hyperparameter $\mathcal{T} \in \mathbb{T}$ and sets $\mathcal{L}(\mathcal{T})$ and $\mathcal{N}(\mathcal{T})$. The last two steps are based on the fact that i) for any (fixed) reconfiguration decision $\mathbf{u}$ that satisfies (1j)-(1l), there always exists a $\mathcal{T} \in \mathbb{T}$ that has a one-to-one correspondence with this $\mathbf{u}$, and ii) reconfiguration decisions for a distribution network are finite.

Based on Lemma 4, we immediately derive

*Corollary 2:* Let (A.1) = (1a). For every locally optimal solution of (1) with load $\mathbf{d}^* \in \mathcal{D}$, i.e., $col(\mathbf{z}^*, \mathbf{u}^*)$, there always exists a $\mathcal{T}^* \in \mathbb{T}$ such that this $\mathbf{z}^*$ is also a locally optimal solution of $\mathbf{P}_{R\_E}(\mathbf{d}^*, \mathcal{T}^*)$.

Next, we analyze the existence of Lagrange multipliers at the optimal solution of DNR model (1).

*Proposition 1:* If the assumption in Lemma 3 holds, then, there exist Lagrange multipliers for DNR model (1) at every locally optimal solution almost surely such that the KKT conditions are satisfied.

*Proof:* Without loss of generality, let (A.1) = (1a), let $col(\mathbf{z}^*, \mathbf{u}^*)$ denote the locally optimal solution of (1) with load $\mathbf{d}^*$ ($\mathbf{d}^* \in \mathcal{D}$), and let $\mathcal{T}^*$ ($\mathcal{T}^* \in \mathbb{T}$) correspond to this $\mathbf{u}^*$. Based on Corollary 2, this $\mathbf{z}^*$ is also a locally optimal solution of $\mathbf{P}_{R\_E}(\mathbf{d}^*, \mathcal{T}^*)$. According to Lemma 3, we know that $\mathbf{P}_{R\_E}(\mathbf{d}^*, \mathcal{T}^*)$ satisfies the LICQ at every feasible solution, including $\mathbf{z}^*$, almost surely. Mathematically, there exist $\mathbf{\iota}_P$ and $\mathbf{\kappa}_P$ such that the KKT conditions

$$-\nabla w(\mathbf{z}^*)^T = \mathbf{\iota}_P{}^T col(\nabla \mathbf{g}_1(\mathbf{z}^*), \nabla \mathbf{g}_2(\mathbf{z}^* \mid \mathbf{d}^*)) + \mathbf{\kappa}_P{}^T \nabla \mathbf{h}_1(\mathbf{z}^*) \quad (8a)$$
$$\mathbf{g}_1(\mathbf{z}^*) = \mathbf{0}, \mathbf{g}_2(\mathbf{z}^* \mid \mathbf{d}^*) = \mathbf{0}, \quad (8b)$$
$$\mathbf{h}_1(\mathbf{z}^*) \leq \mathbf{0}, \quad (8c)$$
$$\mathbf{\kappa}_P \geq \mathbf{0}, \quad (8d)$$
$$\mathbf{\kappa}_P{}^T \mathbf{h}_1(\mathbf{z}^*) = \mathbf{0}. \quad (8e)$$

hold almost surely for $\mathbf{z}^*$, i.e., the locally optimal solution of $\mathbf{P}_{R\_E}(\mathbf{d}^*, \mathcal{T}^*)$. Now, we explore the existence of Lagrange multipliers in (1). For the locally optimal solution $col(\mathbf{z}^*, \mathbf{u}^*)$ of (1), we derive matrix (11). The elements in matrix (11) denote the partial derivatives of constraints (listed in the row above the matrix) with respect to variables (listed in the column on the left-hand side of the matrix), e.g., $\partial(1d)/\partial p_{mn} = 2r_{mn}$. The row below matrix (11) consists of the Lagrange multipliers for corresponding constraints. For DNR model (1), let $\mathbf{\kappa}^d = col(\kappa^d_{mn}), \forall mn \in \mathcal{L}$. By employing this rule, we can derive other Lagrange multipliers, e.g., $\mathbf{\kappa}^e, \mathbf{\kappa}^{g+}$, and $\mathbf{\iota}^j$. Let $\partial(1d) = col(0, 0, 2r_{mn}, 2x_{mn}, -(r^2_{mn} + x^2_{mn}), -1, 1, M)$, i.e., $\partial(1d)$ denotes a vector, whose Lagrange multiplier is $\kappa^d_{mn}$. Similarly, we derive other vectors, e.g., $\partial(1e)$ and $\partial(1g)^+$. For $\mathbf{P}_{R\_E}(\mathbf{d}, \mathcal{T})$, we use $\iota^4_{mn}$ to denote the Lagrange multiplier of line $mn$ in (A.4) (similar to $\kappa^d_{mn}$ for line $mn$ in (1d)), and let $\mathbf{\iota}^4_P = col(\iota^4_{mn}), \forall mn \in \mathcal{L}(\mathcal{T}), m \in \mathcal{N}(\mathcal{T})$. By employing this rule, we derive other Lagrange multipliers, e.g., $\mathbf{\iota}^6_P$ and $\mathbf{\kappa}^{8+}_P$. Let $\partial(A.4) = col(0, 0, 2r_{mn}, 2x_{mn}, -(r^2_{mn} + x^2_{mn}), -1, 1)$, i.e., $\partial(A.4)$ denotes a vector, whose Lagrange multiplier is $\iota^4_{mn}$. Similarly, we derive other vectors, e.g., $\partial(A.6)$ and $\partial(A.8)^+$.

First, we focus on $\mathbf{\kappa}^d$ and $\mathbf{\kappa}^e$. For switched-on lines, e.g., line $mn$ ($u^*_{mn} = 1$), we have

$$\partial(1d) - M\partial(1j) = col(0, 0, 2r_{mn}, 2x_{mn}, -(r^2_{mn} + x^2_{mn}), -1, 1, 0)$$
$$= col(\partial(A.4), 0), \quad (9a)$$
$$\partial(1e) - M\partial(1j) = col(0, 0, -2r_{mn}, -2x_{mn}, r^2_{mn} + x^2_{mn}, 1, -1, 0)$$
$$= -col(\partial(A.4), 0). \quad (9b)$$

Accordingly, we set $\kappa^d_{mn} = \iota^4_{mn}, \kappa^e_{mn} = 0, \iota^{j\_d}_{mn} = -M\iota^4_{mn}, \iota^{j\_e}_{mn} = 0$ or $\kappa^d_{mn} = 0, \kappa^e_{mn} = -\iota^4_{mn}, \iota^{j\_d}_{mn} = 0, \iota^{j\_e}_{mn} = M\iota^4_{mn}$, depending on the (positive or negative) value of $\iota^4_{mn}$. For switched-off lines, e.g., line $ij$ ($u^*_{ij} = 0$), we set $\kappa^d_{ij} = \kappa^e_{ij} = \iota^{j\_d}_{ij} = \iota^{j\_e}_{ij} = 0$. By following this method, we derive $\mathbf{\kappa}^d \geq \mathbf{0}$ and $\mathbf{\kappa}^e \geq \mathbf{0}$.

Second, we focus on $\mathbf{\kappa}^{g+}$ and $\mathbf{\kappa}^{g-}$. For switched-on lines, e.g., line $mn$ ($u^*_{mn} = 1$), we have

$$\partial(1g)^+ + P^{max}_{mn}\partial(1j) = col(0, 0, 1, 0, 0, 0, 0, 0) = col(\partial(A.8)^+, 0), \quad (10a)$$
$$\partial(1g)^- + P^{max}_{mn}\partial(1j) = col(0, 0, -1, 0, 0, 0, 0, 0) = col(\partial(A.8)^-, 0). \quad (10b)$$

| / | $\partial(1d)$ | $\partial(1e)$ | $\partial(1g)$ | | $\partial(1h)$ | | $\partial(1j)$ | |
|---|---|---|---|---|---|---|---|---|
| $\partial p_n$ | 0 | 0 | 0 | 0 | 0 | 0 | 0 | |
| $\partial q_n$ | 0 | 0 | 0 | 0 | 0 | 0 | 0 | |
| $\partial p_{mn}$ | $2r_{mn}$ | $-2r_{mn}$ | 1 | -1 | 0 | 0 | 0 | |
| $\partial q_{mn}$ | $2x_{mn}$ | $-2x_{mn}$ | 0 | 0 | 1 | -1 | 0 | |
| $\partial l_{mn}$ | $-(r^2_{mn} + x^2_{mn})$ | $r^2_{mn} + x^2_{mn}$ | 0 | 0 | 0 | 0 | 0 | (11) |
| $\partial v_m$ | -1 | 1 | 0 | 0 | 0 | 0 | 0 | |
| $\partial v_n$ | 1 | -1 | 0 | 0 | 0 | 0 | 0 | |
| $\partial u_{mn}$ | M | M | $-P^{max}_{mn}$ | $-P^{max}_{mn}$ | $-Q^{max}_{mn}$ | $-Q^{max}_{mn}$ | $2u^*_{mn} - 1$ | |
| | $\kappa^d_{mn}$ | $\kappa^e_{mn}$ | $\kappa^{g+}_{mn}$ | $\kappa^{g-}_{mn}$ | $\kappa^{h+}_{mn}$ | $\kappa^{h-}_{mn}$ | $\iota^j_{mn}$ | |

Accordingly, we set $\kappa_{mn}^{g+} = \kappa_{mn}^{8+}$, $\kappa_{mn}^{g-} = \kappa_{mn}^{8-}$, $\iota_{mn}^{j\_g+} = P_{mn}^{max}\kappa_{mn}^{8+}$, $\iota_{mn}^{j\_g-} = P_{mn}^{max}\kappa_{mn}^{8-}$. For switched-off lines, e.g., line $ij$ ($u_{ij}^* = 0$), we have

$$\partial(1g)^+ + P_{ij}^{max}\partial(1j) = col(0, 0, 1, 0, 0, 0, 0, 0) = col(\partial(A.6), 0), \quad (12a)$$

$$\partial(1g)^- + P_{ij}^{max}\partial(1j) = col(0, 0, -1, 0, 0, 0, 0, 0) = -col(\partial(A.6), 0). \quad (12b)$$

Accordingly, we set $\kappa_{ij}^{g+} = \iota_{ij}^{6}$, $\kappa_{ij}^{g-} = 0$, $\iota_{ij}^{j\_g+} = P_{ij}^{max}\iota_{ij}^{6}$, $\iota_{ij}^{j\_g-} = 0$ or $\kappa_{ij}^{g+} = 0$, $\kappa_{ij}^{g-} = -\iota_{ij}^{6}$, $\iota_{ij}^{j\_g+} = 0$, and $\iota_{ij}^{j\_g-} = -P_{ij}^{max}\iota_{ij}^{6}$, depending on the (positive or negative) value of $\iota_{ij}^{6}$. By following this method, we derive $\boldsymbol{\kappa}^{g+} \geq \mathbf{0}$ and $\boldsymbol{\kappa}^{g-} \geq \mathbf{0}$. Similarly, we derive $\boldsymbol{\kappa}^{h+} \geq \mathbf{0}$ and $\boldsymbol{\kappa}^{h-} \geq \mathbf{0}$.

Third, we focus on $\boldsymbol{\iota}^{j}$. Based on the above analysis, we have $\iota_{mn}^{j} = \iota_{mn}^{j\_d} + \iota_{mn}^{j\_e} + \iota_{mn}^{j\_g+} + \iota_{mn}^{j\_g-} + \iota_{mn}^{j\_h+} + \iota_{mn}^{j\_h-}$, $\forall mn \in \mathcal{L}$. By following this method, we derived $\boldsymbol{\iota}^{j}$.

Last, we focus on the other Lagrange multipliers, i.e., $\boldsymbol{\iota}^{b}$, $\boldsymbol{\iota}^{c}$, $\boldsymbol{\iota}^{f}$, $\boldsymbol{\kappa}^{i+}$, $\boldsymbol{\kappa}^{i-}$, $\boldsymbol{\kappa}^{k}$, and $\boldsymbol{\iota}^{l}$. Since $\partial(1b) = \partial(A.2)$, $\partial(1c) = \partial(A.3)$, $\partial(1f) = \partial(A.5)$, $\partial(1i)^+ = \partial(A.10)^+$, $\partial(1i)^- = \partial(A.10)^-$, we set $\boldsymbol{\iota}^{b} = \boldsymbol{\iota}_{P}^{2}$, $\boldsymbol{\iota}^{c} = \boldsymbol{\iota}_{P}^{3}$, $\boldsymbol{\iota}^{f} = \boldsymbol{\iota}_{P}^{5}$, $\boldsymbol{\kappa}^{i+} = \boldsymbol{\kappa}_{P}^{10+} \geq \mathbf{0}$, $\boldsymbol{\kappa}^{i-} = \boldsymbol{\kappa}_{P}^{10-} \geq \mathbf{0}$. In addition, we set $\boldsymbol{\kappa}^{k} = \mathbf{0}$ and $\boldsymbol{\iota}^{l} = \mathbf{0}$.

Overall, we derive Lagrange multipliers for DNR at locally optimal solution $col(\mathbf{z}^*, \mathbf{u}^*)$ almost surely such that the KKT conditions hold. This completes the proof. ∎

*Remark 4:* Despite that the constraints in (1) at any locally optimal solution do not satisfy the LICQ, according to Proposition 1 and [13], they may satisfy other weaker constraint qualifications, verifying the existence of Lagrange multipliers. In other words, there exist Lagrange multipliers for (1) at locally optimal solution $col(\mathbf{z}^*, \mathbf{u}^*)$ almost surely, and we derive one of them in Proposition 1. More importantly, Proposition 1 establishes a solid theoretical foundation for solving *general* DNR models via KKT conditions-based methods, e.g., Newton's methods [12]. In addition, according to [17], the convergence of the distributed algorithm ALADIN requires the existence of Lagrange multipliers. This proposition also helps analyze the convergence of the distributed algorithm ALADIN and its application for solving distributed DNR problems.

TABLE II
OPTIMALITY OF RECONFIGURATION SOLUTIONS

|  | Topo. 1 | Topo. 2 | Topo. 3 |
|---|---|---|---|
| Reconfig. solution | [1, 1, 1, 0] | [1, 1, 0, 1] | [1, 0, 1, 1] |
| Obj. value (p.u.) | 6.72e-3 | **3.64e-3** | 8.61e-3 |
| Opt. reconfig. solution | Topo. 2 [1, 1, 0, 1] (errors <1e-8) |  |  |
| Opt. obj.value (p.u.) | 3.64e-3 |  |  |

## IV. CASE STUDIES

This section conducts numerical tests on modified 4-bus and 33-bus distribution networks. The small-scale 4-bus network is used to demonstrate the optimality of the reconfiguration solution and the existence of Lagrange multipliers in DNR problems. The medium-scale 33-bus network is used to validate the scalability of DNR model (1), as well as the existence of Lagrange multipliers. The parameters and network topologies of both networks are listed in [18]. The code is written in Julia 1.11.3 with JuMP and implemented on a laptop with an Intel(R) i58265U CPU. The solver for MINCPs and NCPs is SCIP 10.2.0.

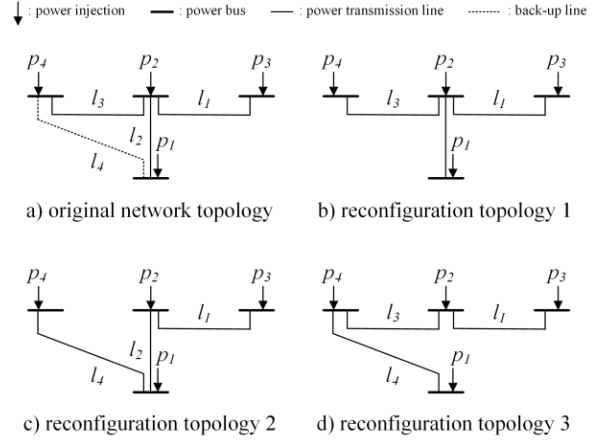

Fig. 2. Four-bus distribution network topology and all possible reconfiguration topologies.

First, we test the optimality of the reconfiguration solution derived by solving model (1). For ease of comparison, we present the network topology of the 4-bus distribution network in Fig. 2 a) and all its possible reconfiguration topologies in Fig. 2 b), c), and d), each corresponding to a reconfiguration solution. Table II lists their (locally) optimal objective values, i.e., Obj. value, and line states, i.e., Reconfig. solution. Since we enumerate all possible network topologies of this small-scale network, the optimal reconfiguration solution derived by solving model (1), i.e., Opt. reconfig. solution, is validated to be the global optimum, i.e., Opt. obj.value, with high accuracy (i.e., 1e-8). The error is reasonable due to numerical issues. Based on this comparison, we also observe that network topology has a significant impact on network losses, rendering system operators an effective approach for achieving economic system operations.

TABLE III
SCALABILITY OF DNR MODEL (1)

| Computation time (s) | 4-bus network | 33-bus network |
|---|---|---|
| Numbering for off state lines | 3 | 3, 9, 14 32, 37 |
| DNR model (1) (NCP) | 1.33e-1 | 68.08 |
| Original DNR model (MINCP) | 3.60e-2 | 151.36 |

Then, we test the existence of Lagrange multipliers for the (globally optimal reconfiguration) solution. For this solution, we derive Lagrange multipliers that satisfy the KKT conditions. Please refer to [18] for detailed data. This test numerically validates the existence of Lagrange multipliers for DNR problems, i.e., Proposition 1.

Next, we validate the scalability of DNR model (1), an NCP, by comparing it with the original DNR model, an MINCP. The original DNR model is derived by removing constraint (1j) in model (1) and replacing variable $\mathbf{z} \in \mathbb{R}^{3|\mathcal{N}|+4|\mathcal{L}|} = col(p_j, q_j, p_{ij}, q_{ij}, l_{ij}, v_j, u_{ij})$, $\forall ij \in \mathcal{L}$, $j \in \mathcal{N}$, in model (1) by $\mathbf{z}' \in \mathbb{R}^{3|\mathcal{N}|+3|\mathcal{L}|} \times \{0, 1\}^{|\mathcal{L}|} = col(p_j, q_j, p_{ij}, q_{ij}, l_{ij}, v_j, u_{ij})$, $\forall ij \in \mathcal{L}$, $j \in \mathcal{N}$, i.e., $u_{ij} \in \{0, 1\}$. Namely, the original DNR model consists of (1a)-(1i), (1k), (1l) with variable $\mathbf{z}'$. Test results are shown in Table III. In this table, computation time is the executive time of JuMP function "JuMP.optimize!" Based on the results, DNR model (1) and the original DNR model derive the same (globally/locally) optimal reconfiguration solutions for 4- and 33-bus distribution networks, respectively. Thus, we do not need to consider the impact of different solutions on computation time. For the small-scale 4-bus distribution network, the original DNR model outperforms DNR model (1) in terms of computation time. This is reasonable as the SCIP solver solves an MINCP by using a branch-and-cut method [19]. A small number of binary line state variables does not need extensive branch-and-cut processes. However, when the number of binary variables increases, e.g., the 33-bus distribution network, the limitation of solving a (complex) MINCP manifests. Based on the results in Table III, it induces a much longer computation time compared to solving an NCP. This case validates the scalability of DNR model (1). At last, we also test the existence of Lagrange multipliers for the optimal solution of DNR model (1) and accordingly derive Lagrange multipliers, validating their existence. Please refer to [18] for detailed data.

## V. CONCLUSION

In this brief, we study the existence of Lagrange multipliers in DNR problems. First, we theoretically prove that the LICQ does not hold for DNR model (1) at any feasible solution, indicating that unique Lagrange multipliers do not exist for this model. Interestingly, we find that it is the reconfiguration constraints that induce the violation the LICQ. Then, we prove that under mild conditions, there exist Lagrange multipliers for DNR model (1) at every locally optimal solution almost surely such that the KKT conditions are satisfied. This conclusion theoretically supports the implementation of various KKT-based algorithms for solving the DNR model, e.g., Newton's method and the interior point method. Case studies demonstrate that i) model (1) has the potential to be solved to a globally optimal reconfiguration solution, ii) network topology has a significant impact on system losses and should be paid great attention, and iii) compared with the original DNR model, i.e., an MINCP, model (1), an NCP, exhibits better scalability of for medium-scale DNR problems. Future work includes studying distributed DNR problems.

## APPENDIX

### A. Formulation of the Auxiliary Model

The formulation of $\mathbf{P}_{R\_E}(\mathbf{d}, \mathcal{T})$ is presented as follows.

$$\mathbf{P}_{R\_E}(\mathbf{d}, \mathcal{T}): \min_{\mathbf{z}} w(\mathbf{z}) \tag{A.1}$$

s.t. $p_n = g_n^P - d_n^P = \sum_{mn \in \mathcal{L}_n^+} (r_{mn} l_{mn} - p_{mn}) + \sum_{no \in \mathcal{L}_n^-} p_{no}, \quad \forall n \in \mathcal{N},$ (A.2)

$q_n = g_n^Q - d_n^Q = \sum_{mn \in \mathcal{L}_n^+} (x_{mn} l_{mn} - q_{mn}) + \sum_{no \in \mathcal{L}_n^-} q_{no}, \quad \forall n \in \mathcal{N},$ (A.3)

$v_m - v_n = 2(r_{mn} p_{mn} + x_{mn} q_{mn}) - (r_{mn}^2 + x_{mn}^2) l_{mn},$
$\quad \forall mn \in \mathcal{L}(\mathcal{T}),$ (A.4)

$v_m l_{mn} = p_{mn}^2 + q_{mn}^2, \quad \forall mn \in \mathcal{L},$ (A.5)

$p_{mn} = 0, \quad \forall mn \in \mathcal{L} \backslash \mathcal{L}(\mathcal{T}),$ (A.6)

$q_{mn} = 0, \quad \forall mn \in \mathcal{L} \backslash \mathcal{L}(\mathcal{T}),$ (A.7)

$-P_{mn}^{max} \leq p_{mn} \leq P_{mn}^{max}, \quad \forall mn \in \mathcal{L}(\mathcal{T}),$ (A.8)

$-Q_{mn}^{max} \leq q_{mn} \leq Q_{mn}^{max}, \quad \forall mn \in \mathcal{L}(\mathcal{T}),$ (A.9)

$(V_n^{min})^2 \leq v_n \leq (V_n^{max})^2, \quad \forall n \in \mathcal{N}.$ (A.10)

Objective function $w(\cdot)$ can be any physically practical (1-time) continuous differentiable function, .e.g., objective function (1a) for minimizing aggregated line losses. Hyperparameter $\mathbf{d} \in \mathcal{D} = col(d_n^P, d_n^Q)$, $\forall n \in \mathcal{N}$. Set $\mathcal{D}$ consists of all legal values of loads and is assumed to be a nonempty open set with non-zero measure [13]. Hyperparameter $\mathcal{T} \in \mathbb{T}$ denotes the graph of a connected radial power sub-network, whose edges and vertices consist of $mn \in \mathcal{L}(\mathcal{T})$ and $n \in \mathcal{N}$, respectively. $\mathbb{T}$ denotes the set of all graphs of connected radial power sub-networks. In other words, for any $\mathcal{T} \in \mathbb{T}$, the topology of the corresponding power sub-network must be a tree. Note that $\mathbf{P}_{R\_E}(\mathbf{d}, \mathcal{T})$ is not restricted to radial power networks and can also be used for meshed power networks.